\documentclass[11pt]{article}

\usepackage{amssymb}

\usepackage{latexsym}
\usepackage{amssymb}
\usepackage{amsfonts}
\def\no{\nonumber}
\def\be{\begin{equation}}
\def\ee{\end{equation}}
\def\ba{\begin{eqnarray}}
\def\ea{\end{eqnarray}}

\def\btu{\bigtriangleup}
\def\btd{\bigtriangledown}
\def\e1{\epsilon}
\def\AAl{\mathcal{A}_{\lambda}}
\def\A0{\stackrel{\circ}{\AAl}}

\def\vt{\vartheta}
\def\o1{\Omega}
\def\a{\alpha}
\def\lmd{\lambda}

\def\v1{\varphi}

\def\0v{\varphi_{0}}
\def\d1{\delta}
\def\part{\partial}

\def\f1{\frac}
\def\t1{\theta}
\def\s1{\sqrt{\e1}}
\def\b1{\beta (r,\theta )}
\def\bar{\overline}
\def\bs{\begin{eqnarray*}}
\def\es{\end{eqnarray*}}

\begin{document}

\large{\title {An analysis of the two-vortex case in the Chern-Simons Higgs
model\thanks{{\bf dg-ga/9710006}}}
\author{Weiyue Ding\thanks{Institute of Mathematics,
Academia Sinica, Beijing 100080, P. R. China.}\\
\and
J\"{u}rgen Jost\thanks{Max-Planck-Institute for Mathematics in the
Sciences,
Inselstr. 22-26, 04103 Leipzig, Germany.}\\
\and
Jiayu Li\thanks{Institute of Mathematics,
Academia Sinica, Beijing 100080, P. R. China.}\\
\and
Guofang Wang\thanks{Institute of Mathematics,
Academia Sinica, Beijing 100080, P. R. China.}}
\date {June 1997}
\maketitle

\newtheorem{theorem}{Theorem}[section]
\newtheorem{lemma}[theorem]{Lemma}
\newtheorem{corollary}[theorem]{Corollary}
\newtheorem{remark}[theorem]{Remark}
\newtheorem{definition}[theorem]{Definition}
\newtheorem{proposition}[theorem]{Proposition}
\begin {abstract} Extending work of Caffarelli-Yang and Tarantello, we
present a
variational existence proof for two-vortex solutions of the periodic
Chern-Simons Higgs model and analyze the asymptotic behavior of these
solutions
as the parameter coupling the gauge field with the scalar field tends to
$0$.
\end{abstract}

\section{Introduction}
The present note deals with a simplified form of the so-called
anyon model, a classical field theory defined on (2+1) dimensional
Minkowski space where the Lagrangian couples Maxwell and Chern-Simons terms
coming from a gauge field with a scalar field. Multivortex solutions of
this model are of interest for example in superconductivity.
Since the corresponding Euler-Lagrange equations are somewhat complicated
and since at large distances and low energies, the lower order Chern-Simons
terms
dominates the higher order Maxwell term,
Hong-Kim-Pac [HKP] and Jackiw-Weinberg [JW] proposed to study the
simplified
version without the Maxwell term. Furthermore, they found that for a
special
$6^{th}$ order choice of the potential, there may exist
 time-independent vortex solutions
satisfying a first order Bogomolny type self-dual equations, similar to the
Abelian Higgs equations that arise from a theory with a $4^{th}$ order
potential. Nevertheless, the former theory is mathematically and
physically richer, because the potential admits both symmetric and
asymmetric
minima, and varying
the coupling parameter allows to interpolate between the different types of
vacua.

Caffarelli-Yang [CY] and Tarantello [T] then obtained existence results for
stationary periodic multivortices. Mathematically, this means that the
stationary theory is studied on a two-dimensional torus. While the
existence
results hold for an arbitrary number of prescribed vortices, some of the
finer
variational aspects and asymptotic results for a coupling parameter tending
to
$0$ could be obtained only for the case of a single vortex. In the present
paper, we take up those studies and analyze the two-vortex case. The
two-vortex
case is more difficult because it turns out to be a borderline case for the
Moser-Trudinger inequality. For more than two vortices, special periodic
situations may be reduced to the one or two-vortex case, but the general
situation  is unclear at present.

{\bf Acknowledgement}:{This paper was written at the Max-Planck-Institute
for
Mathematics in the Sciences in Leipzig.  The first author thanks the
Max-Planck-Institute
for its hospitality and good working conditions. The third author
 was supported by a fellowship of the
Humboldt Foundation. The fourth author was supported by the DFG through the
Leibniz award
of the second author. }

\section{Geometrical description of the model}

The so-called Chern-Simons Higgs functional introduced by Hong-Kim-Pac
[HKP] and Jackiw-Weinberg  [JW] has the Lagrangian density
\be
{\mathcal{L}}(\phi,A) = D_{\a}\phi  \overline{D_{\a} \phi}+ \frac{1}{4} k
\varepsilon^{\alpha \beta \gamma} F_{\alpha \beta} A_{\gamma} - V(|\phi|)
\ee
for a complex scalar field $\phi$ coupled with a Chern-Simons gauge field
$A$
on $2+1$ dimensional Minkowski space $\mathbb{R}^{1,2}$. In geometric
terminology, $\phi$ can be considered as a section of the bundle
$\mathbb{R}^{1,2}\times \mathbb{C}$, and $A=-i A_{\a}dx^{\a}$
$(A_{\a}\in\mathbb{R}$, $x=(x^0,x^1,x^2)\in\mathbb{R}^{1,2}$) is a unitary
connection on this bundle. (We normalize the electric charge to be 1 here.)
The curvature of $A$ is
$$F_A = \frac{-i}{2} F_{\alpha \beta}\, dx^{\alpha} \wedge dx^{\beta}$$
with $F_{\alpha\beta} = \partial_{\alpha}A_{\beta} -
\partial_{\beta}A_{\alpha}$, $  \alpha,  \beta =0, 1, 2$,
and we write for the covariant derivative
$$
D_A\phi=D_{\a}\phi dx^{\a} \qquad\hbox{ with }
D_{\a}\phi=\part\phi-iA_{a}\phi.$$
In (1) $V$ is a potential term, and in the model under consideration,
it is taken as

\be
V(|\phi|)=\frac{1}{k^2}|\phi|^2(1-|\phi|^2)^2.
\ee
Here, the coupling parameter $k$ is the same as in ({1}).
This choice of $V$ will lead us to Bogomolny type
selfduality equations, following the derivation in [JW]. Finally,
$\e1^{\a\beta\gamma}$ is normalized by $\e1^{012}=1$, and of course,
the usual summation convention is in force; and we use the Minkowski metric
to raise and lower indices. Thus, all notations in ({1}) have
been explained.

The Euler-Lagrange equations for the Lagrangian action density
$\mathcal{L}$
are
\be
\begin{array}{rcccr}
\frac{1}{2} k\varepsilon^{\alpha \beta \gamma} F_{\alpha \beta} & = &
j^{\gamma}
& = & i(\phi \overline{D^{{\gamma}} \phi} - \bar{\phi} D^{\gamma} \phi) ,
\no\cr
 D_{\a}D^{\a} \phi & = & - \frac{\partial V(\phi)}{\partial \phi},\cr
\end{array}
\ee
where $j^{\gamma}$ is the conserved matter current density.
[HKP] and [JW] are interested in time independent vortex solutions to
these field equations. For such static configuration, the energy becomes
\be
 E(\phi,A) = \int d^2x(|D_A\phi|^2 -A_0^2|\phi|^2-kA_0F^{12} +V(|\phi|))
\ee
Varying {it w.r.t} $A_0$ yields
$$A_0=-\frac{k}{2}\frac{F^{12}}{|\phi|^2},$$
and $E$ becomes
\be
E(\phi,A) = \int d^2x( |D_A \phi|^2 + \frac{k^2}{4}
\frac{|F|^2}{|\phi|^2} + V(|\phi|)).
\ee
With the choice ({2}) for $V(|\phi|)$, ({5}) may be
rewritten as
\be
\begin{array}{rcccr}
E(\phi,A)&=&\int d^2x ( |(D_1\pm D_2)\phi|^2+  (\frac{k}{|\phi|} F_{12}
 \mp \frac{2}{k}
|\phi| (|\phi|^2 -1))^2 \no \\
&&\pm F_{12}+\mbox{Im }\{\partial_j \varepsilon_{jk} \bar{\phi} D_k
\phi\}).
\end{array}
\ee

$E$ is gauge invariant in the sense that its value is unaffected under
changing $(\phi,A)$ to $(e^{i\vartheta}\phi, A_j+\part_j\vartheta)$ with
real
valued $\vartheta$.

Following Caffarelli-Yang [CY], we wish to study solutions that are
periodic
{w.r.t.} some lattice on $\mathbb{R}^2$. Such solutions can
be interpreted as solutions on a fundamental domain
$$\Omega=\{z^1\tau_1+z^2\tau_2, 0<z^1,z^2<1\}$$
for the lattice, where $\tau_1,\tau_2 \in \mathbb{R}^2$ are the generators
of this
lattice, satisfying  so-called 't Hooft boundary conditions that we now
state.
Let $\vt_1,\vt_2$ be (smooth) real-valued functions defined for
$x=z^2\tau_2$
and $x=\tau_1+z^2\tau_2,$ $0<z^2<1$, or $x=z^1\tau_1$
and $x=\tau_2+z^1\tau_1,$ $0<z^1<1$, respectively. We then require on those
boundary lines of $\Omega$ where $\vt_j$ is defined
\be
 e^{i\vt_j(x+\tau_j)}\phi(x+\tau_j)=e^{i\vt_j(x)}\phi(x),
\ee
\be
(A_k+\part_k\vt_j)(x+\tau_j)=(A_k+\part_k\vt_j)(x).
\ee
Since $\phi$ is a single-valued complex function, its phase change
when traversing $\part\Omega$ has to be $2\pi N$ for some integer $N$.
Using
the resulting constraint on the $\vt_j$ from ({7}) in ({8})
yields
\be
\int_{\Omega}{F_{12}}dx=\int_{\part\Omega}A_kdx^k=2\pi N.
\ee
W.l.o.g., we assume $N\ge 0$ in the sequel, in order to obtain the upper
signs in ({6}) in the sequel.
Since the last term in ({6}) is a divergence term, we obtain
\ba
E(\phi,A) & = & \int d^2x \left\{
| (D_1+iD_2)\phi|^2 +(\frac{k}{2}\frac{F_{12}}{|\phi|}-\frac{1}{k} |\phi|
 (|\phi|^2-1))^2 \right\}\no\\
&+& 2\pi N
\ea
The absolute minima of $E$ therefore satisfy the Bogomolny type
self-dual equations
\be
 D_1\phi+iD_2\phi=0
\ee
\be
 k^2F_{12}=2|\phi|^2(1-|\phi|^2)
\ee
subject to the boundary conditions ({7}) and ({8}).

Geometrically, we consider the quotient of $\mathbb{R}^2$ by the lattice
generated by $\tau_1,\tau_2$, i.e. we identify the boundary portion
$z^k\tau_k$, $0\le z_k\le 1$, of $\Omega$ with $z^k\tau_k+\tau_j$, for
$j\not =k\in \{1,2\}$. The resulting torus is denoted by $\Sigma$.
Because of the boundary conditions, $\phi$ can then be considered as
a section of a complex line bundle over $\Sigma$ of degree $N$.
Equation ({11}) then says that $\phi$ is a holomorphic section
of this line bundle, and as such, it has to have $N$ zeros,
counted with multiplicity (unless $\phi\equiv 0$).

Conversely, it was shown by Caffarelli-Yang [CY] that for some critical
number $k_c$ that satisfies the bound
\be
0<k_c\le \frac12(\frac{|\Sigma|}{\pi N})^{\frac12} \qquad (|\Sigma|:=\hbox{
Area}(\Sigma))
\ee
and $0<k<k_c$, and given $p_1,p_2,\cdots,p_N\in \Sigma$
(not necessarily distinct), there exists such a solution of (11) and (12)
with zeros at $p_i$ with corresponding multiplicities in the case some
of the $p_i$ coincide. Tarantello [T] then showed that in this situation,
there exists a second solution different from the one of [CY]. It is of
physical
and mathematical interest to study the asymptotic behavior of these
solutions
as $k$ tends to $0$. Tarantello showed that the Caffarelli-Yang solution
converges
to one in absolute value away from the vortices $p_1,p_2,\cdots,p_N$. In
other
words, in the limit, we obtain a covariantly constant section of a flat
bundle
over $\Sigma\setminus\{p_1,p_2,\cdots,p_N\}$, whereas the curvature becomes
a
sum of delta functions supported on the vortex set. This is analogous to
the
situation in [HJS] where the Ginzburg-Landau theory with a potential term
of
fourth order was investigated. In case $N=1$, Tarantello showed in contrast
that the second solution converges to $0$ uniformly, and furthermore,
after rescaling, one obtains a solution of the mean field equation
\be
\left\{\begin{array}{ll}
\Delta w_0 = -4\pi N\left ( \frac{e^{u_0+w_0}}{\int_{\Sigma}e^{u_0+w_0}}-
\frac{1}
{|\Sigma|}\right )\no \\
\int_{\Sigma}w_0 = 0
\end{array}\right.
\ee
where $u_0\in H^{1,p}(\Sigma)$ ($1<p<2$) is the (negative) Green function
satisfying
\be
\left\{\begin{array}{ll}
\Delta u_0 = -\frac{4\pi N}{|\Sigma|}+4\pi \sum^N_{j=1}\delta_{p_j}\no \\
\int_{\Sigma}u_0 = 0 \\
\end{array}
\right.
\ee

In order to see the reason for the limitation $N=1$ in the
last result, we need to describe the method of the proof.
The aforementioned existence results were obtained by a sub/supersolution
method, i.e.  a method based on the maximum principle. In addition,
Caffarelli-Yang introduced a variational method that we are going to
describe that could reprove the existence results in case $N=1$ only.

We can see that (11) is
\be
\bar\part\phi-\frac{i}{2}(A_1+iA_2)\phi=0
\ee
with $\bar\part=1/2(\part_1+i\part_2)$. Thus, if we know $\phi$,
$A$ is determined by
\be
(A_1+iA_2)=-2i\bar\part\log\phi.
\ee

In order to obtain $\phi$, we put
\be
v(x)=\log|\phi|^2.
\ee

Because of the gauge invariance, given $v$, we may obtain $\phi$
as
\be
\phi(x)=\exp(\frac12v(x)+i\sum^N_{j=1} \arg(x-p_j)).
\ee
Using ({17}), ({12}) translates into the equation
for $v$
\be
\Delta v=\frac{4}{k^2}e^v(e^v-1)+4\pi\sum^N_{j=1}\delta_{p_j}.
\ee
Writing $v=u_0+u$ with the Green function $u_0$ of ({15}), and
putting $K=e^{u_0}$, $\lambda=\frac{4}{k^2}$, we obtain the equation
\be
\Delta u=\lmd Ke^u(Ke^u-1)+\frac{4\pi N}{|\Sigma|}.
\ee

For finding the second solution, Tarantello [T] used the functional
\be
I_{\lmd}(u)=\frac12\int_{\Sigma}(|Du|^2+\frac{\lmd}{2}(Ke^u-1)^2)+
\frac{4\pi N}{|\Sigma|}\int _{\Sigma}u.
\ee
$I_{\lmd}$ is well-defined on $H^{1,2}(\Sigma)$ as
$K=e^{u_0}\in L^{\infty}(\Sigma)$ and since the Moser-Trudinger inequality
([M], [A]) holds:

For every $\e1>0$, there exists $C(\e1)$ with
\be
\log\int_{\Sigma}e^w\le (\frac{1}{16\pi}+\e1)\int_{\Sigma}|Dw|^2
+\frac{1}{|\Sigma|}\int_{\Sigma}w+C(\e1)
\ee
for all $w\in H^{1,2}(\Sigma)$.

One easily verifies that critical points of $I_{\lmd}$ yield solutions of
(21). (The reader may worry about the $\log$ on the l.h.s. of
({23}) when compared with ({22}), but that is handled
by imposing the normalization $\int_{\Sigma}Ke^u=1$.)

The reason why the variational method could made to work to obtain
a solution of ({21}) via a minimization procedure in [CY] and
to find the asymptotic behavior of the second solution in [T] only
in case $N=1$ stems from the constant $\frac{1}{16\pi}$ in ({23)}).
($\frac{1}{16\pi}$ in fact is the optimal value in ({23)}, see
[DJLW1] and [F]). As one sees from comparing ({22}) and
({23}), the case $N=2$ in that integral leads to a critical or
borderline situation for the Moser-Trudinger inequality. Such borderline
cases occur in many problems of geometric and analytic interest, see
 e.g. the monograph [St] of Struwe, and they typically require
an analysis that is both subtle and interesting. We therefore attack the
variational aspects of the $N=2$ vortex case in the present paper. We prove
the following theorem.

\begin{theorem}
Let $\Sigma$ be a flat torus, N=2. The equations (11) and (12) have
two solutions $(A_k^1 ,\phi_k^1 )$ and $(A_k^2 ,\phi_k^2 )$ for
sufficiently small $k$ which
satisfy that $|\phi_k^1 |\to 1$ almost everywhere as $k\to 0$, and that
$|\phi_k^2 |\to 0$ almost everywhere as $k\to 0$.
\end{theorem}

The first solution resembles the topological solution in the non-compact
setting and the second one resembles the non-topological solution.
The first solution $(A_k^1 ,\phi_k^1 )$ was
obtained first by Caffarelli-Yang [CY],
using a sub/supersolution
method. The problem studied in the present paper
is meaningful and interesting also on Riemann surfaces other than a torus.
Thus, in our companion paper [DJLW2], we study the $N=2$ case on $S^2$
where we can find even three solutions with different asymptotic
behaviors. In [DJLW1], we derive some background results
on the Kazdan-Warner problem on higher genus Riemann surfaces of which the
present problem can be considered as a special case.

\section{Proof of the theorem}

For simplicity we assume that $|\Sigma|=1$. The following
Moser-Trudinger inequality on a Riemann surface was proved in [F] and [DJLW1].

\begin{lemma} Let $\Sigma$ be a compact Riemann surface (with a conformal
metric
for which $|\Sigma|=1$).
There exists a constant $C>0$ such that,
$$
\log\int_{\Sigma}e^u\le \frac{1}{16\pi}\int_{\Sigma}|\btd u|^2
+
\int_{\Sigma}u+C
$$
for any $u\in H^{1,2}(\Sigma)$.
\end{lemma}

\hfill $\square$

We also need the following concentration lemma, which can be seen
as a generalization of one on $S^2$ proved by Chang-Yang [ChY]
using the conformal transformations.

\begin{proposition} (Concentration lemma)
Let $M$ be a compact Riemann surface with volume 1.
Given a sequence of functions $u_j\in H^{1,2}(M)$
with $\int _M e^{u_j}=1$ and $\|\btd u_j\|_{L^2(M)}
+16\pi\int _M u_j \leq C$, then either

{\rm (i)}  there is a constant $C_0>0$ such that
$\int_M |\btd u_j |^2
\leq C_0$ or

{\rm (ii)} a subsequence which is also denoted by $u_j$ concentrates
at a point $p\in M$, i.e., for any $r>0$,
$$
\lim_{j\to \infty}\int _{B_r(p)}e^{u_j}=1.
$$
\end{proposition}

Proof: Suppose that (ii) does not hold. That is, every subsequence of
$u_j$ does not concentrate.  Then for any
$p\in M$ there is $0<r<\f1{1}{16}i_M$
($i_M$ is the injectivity radius of $M$) such that
\be
\lim_{j\to \infty}\int _{B_r(p)}e^{u_j}<\d1 <1.
\ee
At least there exists a subsequence $u_{j_{k}}$ satisfying
the above inequality. For simplicity, we shall often implicitly pass
to subsequences in the sequel without mentioning it explicitly.

Since $M$ is compact, we can see that there is a finite set
$\{~(p_l ,r_l )~|~l=1,2,\cdots , L~\}$ satisfying
$$
\lim_{j\to \infty}\int _{B_{r_l}(p_l)}e^{u_j}<\d1 <1.
$$
and $\bigcup _{l=1}^{L}B_{r_l}(p_l)\supset M$.

We assume that
$$
\lim_{j\to \infty}\int _{B_{r_1}(p_1)}e^{u_j}
=\max \{\lim_{j\to \infty}\int _{B_{r_l}(p_l)}e^{u_j}
|l=1,2,\cdots ,L\}.
$$
So, we have
$$
\lim_{j\to \infty}\int _{B_{r_1}(p_1)}e^{u_j}
\geq \alpha _0>0,
$$
for some positive constant $\alpha_0$.

The improved Moser-Trudinger inequality (see Lemma 3.3 below)
 implies that
$$
\lim_{j\to \infty}\int _{M\setminus B_{2r_1}(p_1)}e^{u_j}
=0,
$$
if $\int_M |\btd u_j |^2$ is not bounded.

We choose a normal coordinate system $(x_1,x_2)$ around $p_1$, and we
may assume further that in $B_{16r_1}(p_1)$
$$
\f1{1}{2}|x-y|\leq dist_M(x,y)\leq 2|x-y|.
$$
where $|x-y|=dist_{R^2}(x,y)$.

We consider the square $P_1=\{|x_i|\leq 4r_1 |i=1,2\}
\subset R^2$, and we have
$$
\lim_{j\to \infty}\int _{\exp (P_1)}e^{u_j}
=1.
$$
We divide $P_1$ into $4^2$ equal subsquares.
If $\int_M |\btd u_j |^2$ is not bounded,
using the improved Moser-Trudinger inequality one gets a square $P_2$
which is a union of at most $3^2$  of the equal subsquares, such
that
$$
\lim_{j\to \infty}\int _{\exp (P_2)}e^{u_j}
=1.
$$

Then we can get a sequence of squares $P_n\to p_0\in M$ such that
$$
\lim_{j\to \infty}\int _{B_r(p_0)}e^{u_j}=1
$$
for any $r>0$, which
contradicts (24). This means that
$\int_M |\btd u_j |^2$ is bounded. Therefore (i) holds, this proves the
proposition.

\hfill $\square$

\begin{lemma}([A], [CL])
Let $\Omega_1$ and $\Omega_2$ be two subsets of $\Sigma$ satisfying
{\rm dist}$(\Omega_1,\Omega_2)\ge \epsilon_0>0$ and $\alpha_0\in (0, 1/2)$.
For any $\epsilon>0$, there exists a constant
$c=c(\epsilon,\epsilon_0,\alpha_0)$ such that
$$\int_{\Sigma}e^u\le c\exp\{\frac 1{32\pi-\epsilon}\|\btd u\|^2+
\int_{\Sigma}u \}$$
holds for all $u\in H^{1,2}(\Sigma)$ satisfying
$$ \frac{\int_{\Omega_1}e^u}{\int_{\Sigma}e^u}\ge \alpha_0 \qquad\hbox{ and }
\qquad
\frac{\int_{\Omega_2}e^u}{\int_{\Sigma}e^u}\ge \alpha_0.$$
\end{lemma}

To prove the theorem, we consider the functionals
\bs
J_{\lmd}^{\pm}(w)&=& \f1{1}{2}\int_{\Sigma}|\btd w|^2
-8\pi\log\int_{\Sigma}
Ke^w\\
&& -8\pi\log (1\mp \sqrt{1-B_{\lmd}(w)})-\f1{8\pi}{1\mp
\sqrt{1-B_{\lmd}(w)}}
\es
in the space
$$
{\cal A}_{\lmd}=\{w\in H^{1,2}(\Sigma )|\int_{\Sigma}w=0,B_{\lmd}(w)\leq
1\}
$$
where
$$
B_{\lmd}(w)=\f1{32\pi}{\lmd}\f1{\int _{\Sigma}K^2e^{2w}}
{(\int _{\Sigma}Ke^w)^2}.
$$

\begin{lemma}
Let ${\cal A}_{\lmd}^0$ be the interior of ${\cal A}_{\lmd}$.
If $w^{\pm}_{\lmd}\in{\cal A}_{\lmd}^0$
is a critical point of $J^{\pm}_{\lmd}$,
then
$u^{\pm}_{\lmd}=w^{\pm}_{\lmd}-\log (\f1{\lmd}{16\pi}\int_{\Sigma}
Ke^{w^{\pm}_{\lmd}})
-\log (1\mp \sqrt{1-B_{\lmd}(w^{\pm}_{\lmd})})$
is a solution of (21).
\end{lemma}

Proof:
If $w^{\pm}_{\lmd}\in {\cal A}_{\lmd}^0$ is a critical point of
$J^{\pm}_{\lmd}$
in ${\cal A}_{\lmd}$, then it satisfies
\bs
\btu w^{\pm}_{\lmd}&=&-8\pi (1+\f1{B_{\lmd}(w^{\pm}_{\lmd})}{1\mp
\sqrt{1-B_{\lmd}(w^{\pm}_{\lmd})}})
\f1{Ke^{w^{\pm}_{\lmd}}}{\int _{\Sigma}Ke^{w^{\pm}_{\lmd}}}\\
&&+8\pi\f1{B_{\lmd}(w^{\pm}_{\lmd})}{1\mp
\sqrt{1-B_{\lmd}(w^{\pm}_{\lmd})}}\f1{K^2e^{2w^{\pm}_{\lmd}}}{\int
_{\Sigma}K^2e^{2w^{\pm}_{\lmd}}}+\mu
\es
Integrating yields that $\mu =8\pi$. Computing directly shows that
$u^{\pm}_{\lmd}$ is a solution of (21).
This proves the lemma.

\hfill $\square$

Note that, for any fixed $\lmd >0$, we have
$B_{\lmd}(w)\geq \f1{32\pi}{\lmd}$ by the H\"{o}lder inequality.
We can see that $J^{\pm}_{\lmd}$ are bounded from below
in ${\cal A}_{\lmd}$. Therefore we can get two sequences of functions
$w^{\pm }_{i}\in {\cal A}_{\lmd}$
such that
$$
J_{\lmd}^{\pm}(w^{\pm}_i)\to \inf_{w\in {\cal A}_{\lmd}}
J^{\pm}_{\lmd}(w)
$$
as $i\to\infty$.

We claim that
$w^{\pm}_i$ are bounded in $H^{1,2}(\Sigma)$.
We show this claim for $w^{+}_i$, the proof for $w^{-}_i$ is similar.
It suffices to show that $\|\btd w^{+}_i\|_{L^2(\Sigma)}$ is bounded.

We set $u^{+}_{i}=w^{+}_i+c^{+}_{i}$,  where
$c^{+}_{i}$  is a constant such that
$\int_{\Sigma}e^{u^{+}_i}
=1$.
The boundedness of $J^+_{\lmd}(w^+_i)$ implies
\be\frac12\int_{\Sigma}|\btd u^+_i|^2+8\pi c^+_i\le C\ee
for some constant $C$. By the definition of $c^+_i$, the previous inequality
implies that $u^+_i$ satisfies the conditions in Proposition 3.2. On the other
hand, since $w^{+}_i\in {\cal A}_{\lmd}$, we have
$$
\f1{32\pi}{\lmd}\f1{\int _{\Sigma}K^2e^{2u^{+}_{i}}}
{(\int _{\Sigma}Ke^{u^{+}_{i}})^2}\leq 1.
$$
So
$$
\int _{\Sigma}K^2e^{2u^{+i}_{\lmd}}\leq C_{\lmd}
$$
where $C_{\lmd}$ is a positive constant depending only on
$\lmd$. By the H\"older inequality one gets
$$
\int_{B_{\e1}(p)}Ke^{u^{+}_{i}}\leq C_{\lmd}\e1.
$$
Now applying Proposition 3.2 to the sequence $u^+_i$, we get
the boundedness of $\|\btd u^+_i\|_{L^2(\Sigma)}$, which is equivalent
to our claim.

{}From the claim, we have

(1) $u^+_i$ converges to $u^+_{\lmd}$ weakly in $H^{1,2}(\Sigma)$,

(2) $u^+_i$ converges to $u^+_{\lmd}$ strongly in $L^q(\Sigma)$,
for any $1<q<\infty$,

(3) $\int_{\Sigma}(Ke^{u^+_i})^l\to \int_{\Sigma}(Ke^{u^+_{\lmd}})^l$ as $i\to
\infty$, for $l=1,2,$ by the Moser-Trudinger inequality and the Lebesgue
convergence theorem.

(2) and (3) imply that $u^+_{\lmd}\in {\cal A}_{\lmd}$.
The semi-continuity of the Dirichlet integral, together with (1)--(3), implies
that
$$J^+_{\lmd}(u^+_{\lmd})=\inf_{u\in {\cal A}_{\lmd}}J^+_{\lmd}(u).$$

Therefore, to show the theorem, it suffices to prove that
$J^{\pm}_{\lmd}$ achieves its minimum in ${\cal A}_{\lmd}^0$ if
$\lmd$ is sufficiently large.

By Lemma 3.1 we have
\be
\inf_{w\in \partial{\cal A}_{\lmd}}J^+_{\lmd}(w)\geq -C
-8\pi (1+\log \max K),
\ee
however, choosing $w_0 = 0$, we have
\ba
\inf_{w\in {\cal A}_{\lmd}}J^+_{\lmd}(w)&\leq &
-8\pi\log\int_{\Sigma}K
-8\pi\log (1-\sqrt{1-B_{\lmd}(0)})\no\\
&&-\f1{8\pi}{1-\sqrt{1-B_{\lmd}(0)}}.
\ea
By (26) and (27) we can see that
$J^{+}_{\lmd}$ achieves its minimum in ${\cal A}_{\lmd}^0$ if
$\lmd$ is sufficiently large.





For $J^-_{\lmd}$, with the help of Lemma 3.1 we can show that
\ba
\inf_{w\in \partial {\cal
A}_{\lmd}}J^-_{\lmd}(w)\ge\alpha_0-8\pi+8\pi\log(16\pi)
\ea
and
\ba
\lim_{\lmd\to \infty}\inf_{w\in  {\cal
A}_{\lmd}}J^-_{\lmd}(w)=\alpha_0-4\pi+8\pi\log(8\pi),
\ea
where $\alpha_0=\inf_{u\in H^{1,2}(\Sigma)}\{\frac 12\int_{\Sigma}
|\btd u|^2+8\pi\int_{\Sigma}u-8\pi\log\int_{\Sigma}e^u\}.$
The argument for  (28) and (29) is the same as one given in [DJLW2].
We omit it here. Since $-8\pi+8\pi\log(16\pi)>-4\pi+8\pi\log(8\pi)$,
(28) and (29) imply that
 $J^+_{\lmd}$ achieves its minimum in ${\cal A}^0_{\lmd}$ provided that
 $\lmd$ is sufficiently large.

By Lemma 3.3, we therefore have two solutions of the equations (11) and (12),
$(A^1_k,\phi^1_k)$
and $(A^2_k,\phi^2_k)$ corresponding to $u^{\pm}_{\lmd}$. It was proved
in [T] that $|\phi_k^1|\to 1$ almost everywhere as $k\to 0$.

It is clear that
$$
|\phi^2_k|^2=e^{u_0+u_{\lmd}^-},
$$
so
\bs
\int_{\Sigma}|\phi^2_k|^2&=&\f1{16\pi}{\lmd}\f1{1}
{1+\sqrt{1-B_{\lmd}(w^-_{\lmd})}}\\
&&\leq \f1{16\pi}{\lmd}.
\es
Thus, as $\lambda=\frac 4{k^2}$,
$$
\lim_{k\to 0}\int_{\Sigma}|\phi^2_k|^2=0.
$$
This proves our theorem.

\vspace{.2in}
\begin{center}
{\large\bf REFERENCES}
\end{center}
\footnotesize
\begin{description}

 \item[{[A]}]{Aubin, T., Nonlinear analysis on manifolds, Springer-Verlag,
 New York, 1982.}
 \item[{[BM]}] {Brezis, H. and Merle, F., Uniform estimates and blow up
 behavior for solutions of $-\btu u = V(x)e^u $ in two dimensions,
 Comm. Partial Diff. Equat. 16(1991), 1223-1253.}

 \item[{[CY]}] {Caffarelli, L. and Yang, Y.S., Vortex condensation in the
 Chern-Simons Higgs model: an existence theorem, Comm. Math. Phys.
168(1995),
 321-336.}

 \item[{[ChY]}] {Chang, A. S. Y. and Yang, P., Conformal deformation
 of metrics on $S^2$, J. Diff. Geom. 23(1988), 259-296.}

\item[{[CL]}] {Chen, W. X. and Li, C., Prescribing Gaussian
 curvature on surfaces with conical singularities, J. Geom.
 Anal. 1(1991), 359-372;}

\item[{[DJLW1]}] {Ding, W., Jost, J., Li, J. and Wang, G., The differential
equation $\btu u = 8\pi - 8\pi he^u$ on a
compact Riemann surface, to appear in Asian J. Math. .}

\item[{[DJLW2]}] {Ding, W., Jost, J., Li, J. and Wang, G., Multiplicity
results for the two-vortex Chern-Simons-Higgs model on the
two-sphere, preprint (1997).}
\item[{[F]}] {Fontana, L., Sharp borderline Sobolev inequalities
on compact Riemannian manifolds, Comment. Math. Helv. 68(1993),  415-454.}
\item[{[HJS]}] {Hong, M. C., Jost J. and Struwe,
Asymptotic limits of a Ginzburg-Landau type functional, in: Geometric
Analysis and the
 Calculus of Variations for Stefan Hildebrandt (J. Jost, ed.), Intern.
Press
 Boston, 1996,99-123.}

 \item[{[HKP]}] {Hong, J., Kim, Y. and Pac, P.Y., Multivortex solutions
 of the Abelian Chern-Simons theory, Phys. Rev. Lett. 64(1990), 2230-2233.}

 \item[{[JW]}] {Jackiw, R. and Weinberg, E., Self-dual Chern-Simons
vortices,
 Phys. Rev. Lett. 64(1990), 2234-2237.}

 \item[{[M]}] {Moser, J., A sharp form of an inequality of N. Trudinger,
 Indiana Univ. Math. J. 20(1971), 1077-1092.}


 \item[{[St]}] {Struwe, M. Variational methods, Springer, second edition,
1996.}
 \item[{[T]}] {Tarantello, G., Multiple condensate solutions for the
 Chern-Simons-Higgs theory, J. Math. Phys. 37(1996), 3769-3796.}
 \end{description}

\end{document}